\documentclass[11pt, a4paper]{amsart}
\usepackage{caption}
\usepackage{graphicx}
\usepackage[utf8]{inputenc}
\usepackage{amsmath}
\usepackage{amssymb}
\usepackage{hyperref}
\usepackage{amsthm}
\usepackage{enumerate}
\usepackage{wasysym}

\usepackage{extarrows}
\usepackage{polynom}
\usepackage{dsfont}
\usepackage{verbatim}
\usepackage[shortlabels]{enumitem}
\usepackage[toc,page]{appendix}
\usepackage{a4wide}
\usepackage[english]{babel}
\usepackage{todonotes}
\usepackage{pst-node}
\usepackage{ esint }
\usetikzlibrary{math}
\usetikzlibrary{decorations.pathreplacing}
\theoremstyle{definition}
\newtheorem{theorem}{Theorem}[section]

\newtheorem{definition}[theorem]{Definition}

\newtheorem{remark}[theorem]{Remark}

\newtheorem{conj}[theorem]{Conjecture}
\numberwithin{equation}{section}
\usepackage{tikz-3dplot}
\newcommand{\abs}[1]{\left\lvert#1\right\rvert}

\newcommand{\norm}[1]{\left\|#1\right\|}

\newcommand{\R}{\mathbb R}
\newcommand{\N}{\mathbb N}

\renewcommand{\epsilon}{\varepsilon}

\newcommand{\A}{\mathcal A}
\newcommand{\dx}{\, \mathrm d}
\newcommand{\id}{\mathrm{Id}}

\newcommand{\B}{\mathcal B}

\renewcommand{\phi}{\varphi}
\newcommand{\rev}[1]{\textcolor{black}{#1}}
\usetikzlibrary{hobby}

\usepackage{float}
\begin{document}
\allowdisplaybreaks

\title[On a kinetic Poincar\'e inequality and beyond]{On a kinetic Poincar\'e inequality and beyond}

\date{\today}

\author{Lukas Niebel}
\address[Lukas Niebel]{Institut f\"ur Analysis und Numerik, Universit\"at M\"unster\\
Orl\'eans-Ring 10, 48149 M\"unster, Germany.}
\email{lukas.niebel@uni-muenster.de}
\urladdr{https://lukasniebel.github.io}

\author{Rico Zacher}
\address[Rico Zacher]{Institut f\"ur Angewandte Analysis, Universit\"at Ulm, Helmholtzstra\ss{}e 18, 89081 Ulm, Germany.}
\email{rico.zacher@uni-ulm.de}
\urladdr{https://www.uni-ulm.de/mawi/iaa/members/zacher}

\maketitle

\selectlanguage{english}
\begin{abstract}
{In this article, we give a trajectorial proof of a kinetic Poincar\'e inequality which plays an important role in the De Giorgi-Nash-Moser theory for kinetic equations. The present work improves a result due to J. Guerand and C. Mouhot \cite{guerand_quant_2021} in several directions. We use kinetic trajectories along the vector fields $\partial_t + v \cdot \nabla_x$ and $\partial_{v_i}$, $i = 1,\dots, d$ and do not rely on higher-order commutators such as $[\partial_{v_i},\partial_t + v \cdot \nabla_x] = \partial_{x_i}$ or on the fundamental solution. The presented method also applies to more general hypoelliptic equations. We illustrate this by \rev{investigating} a Kolmogorov equation with $k$ steps.}
\end{abstract}

\vspace{1em}
{\centering \textbf{AMS subject classification.} 35K70, 26D10, 35B45, 35B65, 35H10. \par}
\vspace{1em}

\section{Introduction}

Let $T>0$ and $\Omega_x, \Omega_v \subset \R^d$ be open and bounded sets, we denote $\Omega_T := (0,T) \times \Omega_x \times \Omega_v$. We are concerned with the regularity theory of weak solutions $u= u(t,x,v)$ to the Kolmogorov equation with rough coefficients
\begin{equation} \label{eq:kol}
	(\partial_t + v \cdot \nabla_x) u = \nabla_v \cdot (A \nabla_v u)
\end{equation}
in $\Omega_T$. Here, $A \in L^\infty(\Omega_T ; \R^{d \times d})$ is subject to the following assumptions:
\begin{enumerate}
	\item[(H1)] $\lambda \abs{\xi}^2 \le \langle A(t,x,v) \xi, \xi \rangle$ for all $\xi \in \R^d$ and almost all $(t,x,v) \in \Omega_T$,
	\item[(H2)] $\sum\limits_{i,j = 1}^d \abs{A_{ij}(t,x,v)}^2 \le \Lambda^2$ for almost all $(t,x,v) \in \Omega_T$,
\end{enumerate}
for some constants $0<\lambda <\Lambda$.

	We say that a function $u \in C([0,T];L^2(\Omega_x \times \Omega_v)) \cap L^2((0,T);H_v^1(\Omega_x \times \Omega_v))$ is a weak (super-, sub-) solution to the Kolmogorov equation \eqref{eq:kol} if for all $\varphi \in C^\infty_c(\Omega_T)$ with $\varphi \ge 0$ we have
	\begin{equation} \label{eq:weaksol}
		\int_{\Omega_T} -u (\partial_t + v \cdot \nabla_x ) \varphi+  \langle A \nabla_v u, \nabla_v \varphi \rangle \dx (t,x,v) =  \, (\ge, \, \le) \, 0.
	\end{equation}

Let us give a short overview of the regularity theory related to \eqref{eq:kol}. \rev{Concerning the existence of weak solutions we refer to \cite{ain_weak_2024} and the references therein.} Regularity of weak solutions to equation \eqref{eq:kol} was studied first by Pascucci and Polidoro in \cite{pascucci_moser_2004}, where they obtained local boundedness of weak subsolutions and related mean value inequalities involving powers of $u$. Their work is based on the Moser iteration scheme, \cite{moser_harnack_1964}. H\"older continuity of weak solutions has been obtained first in \cite{wang_nonhom_2007,wang_hypo_2008} by a constructive but complicated method. Recently, a Harnack inequality (which in turn implies H\"older regularity) was proven by Golse, Imbert, Mouhot, and Vasseur in \cite{golse_harnack_2016} by a nonconstructive method based on the ideas of De Giorgi. Later, a first constructive proof of a weak Harnack inequality (which also yields H\"older regularity) was proposed by Guerand and Imbert in \cite{guerand_log_2021}. In \cite{guerand_log_2021}, the authors combine a weak Poincar\'e inequality with the expansion of positivity and an Ink-Spots Theorem to obtain the weak Harnack inequality for some exponent $p>0$. Their method goes back to the work of Moser and Kruzkhov. Local boundedness, together with the weak Harnack inequality, implies the Harnack inequality. The article \cite{guerand_log_2021} gives a very detailed exposition of the historical developments regarding this topic. Guerand and Mouhot have given a different, shorter proof of the weak Harnack inequality in \cite{guerand_quant_2021}. The new key ingredient is a (weak) hypoelliptic Poincar\'e inequality proven by a nice trajectorial argument, which we will explain below. In \cite{guerand_quant_2021}, they follow the ideas of De Giorgi, but all arguments are quantitative in contrast to \cite{golse_harnack_2016}. 

By now there are a lot of results on the Harnack estimate for kinetic equations. However, a precise study which takes into account the exact dependence of the Harnack constant on $\lambda,\Lambda$ is still missing. The authors want to combine $L^p$-estimates, as in \cite{niebel_kinetic_nodate-1}, with a precise Harnack estimate to construct (global) $L^p$-solutions to nonlinear kinetic equations in the future. 
 
Let us now state the (weak) hypoelliptic Poincar\'e inequality of \cite[Proposition 13]{guerand_quant_2021}. The result also allows for a source term and lower order terms, which we drop here for the sake of simplicity.  

\begin{theorem} \label{thm:poincare0}
	Assume (H1) and (H2). Let $\epsilon \in (0,1)$, $\sigma \in (0,1/3)$. Then any nonnegative weak subsolution $u$ to \eqref{eq:kol} in $Q_5 = (-25,0] \times B_{125}(0) \times B_5(0)$ satisfies
	\begin{equation*}
		\norm{(u-\langle u\rangle_{Q_1^-})_+}_{L^1(Q_1)} \le C\left( \epsilon^{-d-2} \norm{\nabla_v u}_{L^1(Q_5)} + \epsilon^\sigma \left(\frac13-\sigma\right)^{-1} \norm{u}_{L^2(Q_5)} \right),
	\end{equation*}
	where $Q_1^- = (-3,-2] \times B_1(0) \times B_1(0)$, $Q_1 = (-1,0] \times B_1(0) \times B_1(0)$, $\langle u\rangle_{Q_1^-} = \fint_{Q_1^-} u \dx(t,x,v)  = \frac{1}{\abs{Q_1^-}} \int_{Q_1^-} u \dx (t,x,v)$ and $C = C(d,\lambda,\Lambda)$.
\end{theorem}

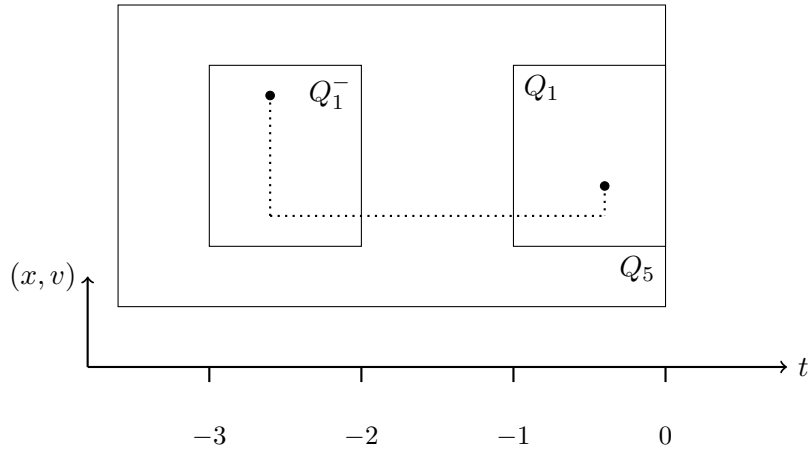
\begin{figure}[H]
\centering
\tikzmath{\x = 0.2; \y = 0.2;}
	\begin{tikzpicture}[scale=4]
  \draw[thick,->] (0, 0) -- (2.3,0) node[right] {$t$};
  \draw[thick,->] (0, 0) -- (0,0.3) node[left] {$(x,v)$};
  \draw[draw=black] (0.3-\x,0.4-\y) rectangle ++(1.8,1);
  \draw[draw=black] (1.6-\x,0.6-\y) rectangle ++(0.5,0.6);
  \draw (1.6-\x,1.2-\y) node[anchor=north west] {$Q_1$};
  \draw[draw=black] (0.6-\x,0.6-\y) rectangle ++(0.5,0.6);
   \draw (1.1-\x,1.2-\y) node[anchor=north east] {$Q_1^-$};
      \draw (2.1-\x,0.6-\y) node[anchor=north east] {$Q_5$};
      \draw [thick] (1.1-\x, 0) -- ++(0, -.05) ++(0, -.15) node [below, outer sep=0pt, inner sep=0pt] {\small\(-2\)};
      \draw [thick] (0.6-\x, 0) -- ++(0, -.05) ++(0, -.15) node [below, outer sep=0pt, inner sep=0pt] {\small\(-3\)};
      \draw [thick] (1.6-\x, 0) -- ++(0, -.05) ++(0, -.15) node [below, outer sep=0pt, inner sep=0pt] {\small\(-1\)};
      \draw [thick] (2.1-\x, 0) -- ++(0, -.05) ++(0, -.15) node [below, outer sep=0pt, inner sep=0pt] {\small\(0\)};
      \filldraw (0.8-\x,1.1-\y) circle[radius=0.4pt];
      \filldraw (1.9-\x,0.8-\y) circle[radius=0.4pt];
      \draw [dotted,thick,  domain=0:1, samples=40] plot ({1.7},{0.6+\x*((0.55+0.9*0.001)/1.1-0.6)});
      \draw [dotted,thick,  domain=0:1, samples=40] plot ({0.6+\x*(1.7-0.6)},{(0.55+0.9*0.001)/1.1});
      \draw [dotted,thick,  domain=0:1, samples=40] plot ({0.6},{(0.55+0.9*0.001)/1.1 + \x*(0.9-(0.55+0.9*0.001)/1.1)});
 \end{tikzpicture}
 \caption{The sets $Q_{5},Q_1,Q_1^-$ in Theorem \ref{thm:poincare0} and the depiction of the projection on the $v$-variable of a trajectory used in the proof of the theorem.}
\end{figure}

Note that employing the underlying scaling $r \mapsto (r^2t,r^3x,rv)$ and translation $(t_0,x_0,v_0) \mapsto (t-t_0,x-x_0-(t-t_0)v_0,v-v_0)$ of the Kolmogorov equation this statement can be transferred to any kinetic cylinder. 

\rev{The proof of Guerand and Mouhot is based on a trajectorial argument, where two points $(t,x,v)$ and $(s,y,w)$ are connected by moving along the vector fields $\partial_t + v \cdot \nabla_x$, $\nabla_v$ and $\nabla_x$. Observe that any two points can always be joined by the trajectory
\begin{equation*}
	(t,x,v) \underset{\nabla_v}{\longrightarrow} \left(t,x, \frac{y-x}{s-t}\right) \underset{\partial_t + v \cdot \nabla_x}{\longrightarrow} \left(s,y, \frac{y-x}{s-t}\right) \underset{\nabla_v}{\longrightarrow} (s,y,w).
\end{equation*}
This leads, however, to problems with integrability when changing variables. They overcome this by adding a subtrajectory of length $\epsilon$ along the vector field $\nabla_x$, i.e. }
\begin{equation*}
	(t,x,v) \underset{\nabla_x}{\longrightarrow} (t,x+\epsilon w,v) \underset{\nabla_v}{\longrightarrow} \left(t,x+\epsilon w, \frac{x+\epsilon w-y}{t-s}\right) \underset{\partial_t + v \cdot \nabla_x}{\longrightarrow} \left(s,y,\frac{x+\epsilon w-y}{t-s}\right) \underset{\nabla_v}{\to} (s,y,w).
\end{equation*}
The trajectory along $\partial_t + v\cdot \nabla_x$ can be controlled using the subsolution property of $u$. Moreover, the trajectory along $\nabla_v$ is controlled by the gradient term on the right-hand side. To control the trajectory along the $\nabla_x$ vector field, an additional result on regularity transfer, compare \cite[Proposition 11]{guerand_quant_2021}, is needed. We stress that assumption (H1) is crucial for their argument to work, see \cite[Proposition 9 and Proposition 11]{guerand_quant_2021}.

In this note, we prove the following kinetic Poincar\'e inequality. 

\begin{theorem} \label{thm:poincare}
	Assume (H2). There exists a universal constant $R_0 \in [1,7]$ such that the following holds. Let $\kappa >0$, $\epsilon \in (0,1)$ and $u$ be a nonnegative subsolution to equation \eqref{eq:kol} in $\tilde{Q} = (-2-\kappa,0] \times B_{R_0(1+\kappa)} \times B_{R_0(1+\kappa^{-1})}$. Then 
	\begin{equation*}
		\norm{(u-\langle u \rangle_{Q_1^-})_+}_{L^1(Q_1)} \le C\left( \epsilon^{-1} \norm{\nabla_v u}_{L^1(\tilde{Q})} + \epsilon^{2d} \norm{u}_{L^1(Q_1)}\right),
	\end{equation*}
	where $Q_1^- = (-2-\kappa,-1-\kappa] \times B_{1} \times B_{1}$, $Q_1 = (-1,0] \times B_1 \times B_1$, $\langle u\rangle_{Q_1^-} = \fint_{Q_1^-} u \dx(t,x,v)  = \frac{1}{\abs{Q_1^-}} \int_{Q_1^-} u \dx (t,x,v)$ and $C = C(d,\kappa,\Lambda)$.
\end{theorem}

Again the result is stated on unit size. By scaling and transformation, we can state it for arbitrary kinetic cylinders. 

For our proof, we use trajectories, which move along the vector fields $\partial_t +v \cdot \nabla_x$ and $\partial_{v_i}$, $i = 1, \dots, d$ only. By that we mean a function $\gamma = (\gamma_t,\gamma_x,\gamma_v) \colon [0,1] \to \R^{1+2d}$ connecting two points $\gamma(0) = (s,y,w)$ and $\gamma(1) = (t,x,v)$ with 
\begin{equation*}
	\frac{\dx}{\dx r} u(\gamma(r)) = \dot{\gamma}_t(r)[(\partial_t + v \cdot \nabla_x) u](\gamma(r)) + \dot{\gamma}_v(r) \cdot [\nabla_v u](\gamma(r))
\end{equation*}
for sufficiently smooth $u \colon \R^{1+2d} \to \R$. Equivalently, $\dot{\gamma}_x(r) = \dot{\gamma}_t(r) \gamma_v(r)$. \rev{We call $\gamma$ a \textit{kinetic trajectory}.}

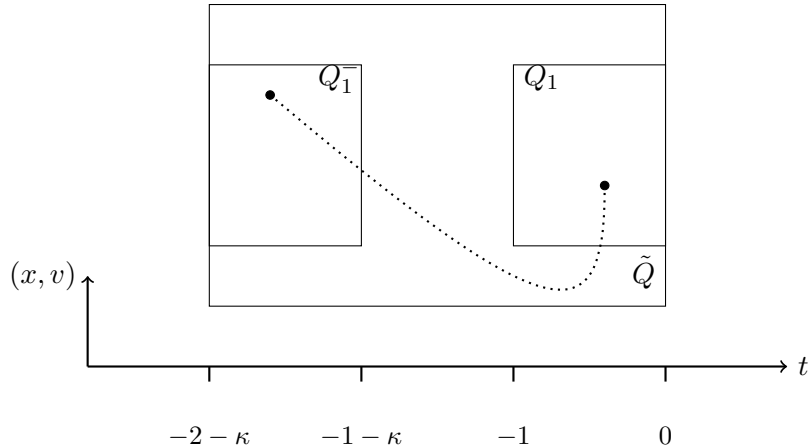
\begin{figure}[H]
\centering
\tikzmath{\x = 0.2; \y = 0.2;}
	\begin{tikzpicture}[scale=4]
  \draw[thick,->] (0, 0) -- (2.3,0) node[right] {$t$};
  \draw[thick,->] (0, 0) -- (0,0.3) node[left] {$(x,v)$};
  \draw[draw=black] (0.6-\x,0.4-\y) rectangle ++(1.5,1);
  \draw[draw=black] (1.6-\x,0.6-\y) rectangle ++(0.5,0.6);
  \draw (1.6-\x,1.23-\y) node[anchor=north west] {$Q_1$ \vphantom{$Q_1^-$}};
  \draw[draw=black] (0.6-\x,0.6-\y) rectangle ++(0.5,0.6);
   \draw (1.13-\x,1.23-\y) node[anchor=north east] {$Q_1^-$};
      \draw (2.1-\x,0.6-\y) node[anchor=north east] {$\tilde{Q}$};
      \draw [thick] (1.1-\x, 0) -- ++(0, -.05) ++(0, -.15) node [below, outer sep=0pt, inner sep=0pt] {\small\(-1-\kappa\)};
      \draw [thick] (0.6-\x, 0) -- ++(0, -.05) ++(0, -.15) node [below, outer sep=0pt, inner sep=0pt] {\small\(-2-\kappa\)};
      \draw [thick] (1.6-\x, 0) -- ++(0, -.05) ++(0, -.15) node [below, outer sep=0pt, inner sep=0pt] {\small\(-1\)};
      \draw [thick] (2.1-\x, 0) -- ++(0, -.05) ++(0, -.15) node [below, outer sep=0pt, inner sep=0pt] {\small\(0\)};
      \filldraw (0.8-\x,1.1-\y) circle[radius=0.4pt];
      \filldraw (1.9-\x,0.8-\y) circle[radius=0.4pt];
      \draw [dotted,thick,  domain=0:1, samples=40] plot ({1.7+\x^2*(0.6-1.7)},{-21/2.2*(\x-\x^(3/2))*0.6+(7*\x^(3/2)-6*\x)*0.9+21/2.2*(\x-\x^(3/2))*1.15+(1-9*\x/2+7/2*\x^(3/2))*0.6});
 \end{tikzpicture}
 \caption{The sets $\tilde{Q},Q_1,Q_1^-$ in Theorem \ref{thm:poincare} and the depiction of the projection on the $v$-variable of a trajectory used in the proof of the theorem.}
\end{figure}

As our kinetic trajectories do not move along the vector field $\nabla_x$, we do not need to use any result on regularity transfer. Hence, there is no additional parameter $\sigma$, and we do not need to impose (H1). More precisely, the matrix does not even need to be positive semidefinite. Moreover, our error term is only in $L^1(Q_1)$ instead of $L^2(\tilde{Q})$ and the gradient term grows like $\epsilon^{-1}$ as $\epsilon \to 0$ only. \rev{The parameter $\kappa$ measures the distance between the two cylinders $Q_1^-$ and $Q_1$.} We note that $\kappa = 1$ leads to the cylinders considered in Theorem \ref{thm:poincare0}. Finally, we emphasise that we need control of the gradient only on $\tilde{Q}$, i.e., no control is needed at times before $Q_1^-$.

The dependence of the size of $\tilde{Q}$ on $\kappa$ is very natural. The closer in time two points are, the more flexibility in velocity we need to connect them while respecting the vector fields $\partial_t +v \cdot \nabla_x$ and $\partial_{v_i}$, $i = 1, \dots, d$. We note that $R_0 $ can be calculated explicitly and turns out to be $R_0 \le 7$ \rev{with our construction}. 

As in \cite{guerand_quant_2021}, we could consider lower order terms and a source term in equation \eqref{eq:kol}. Moreover, one could extend the statement to $L^p$-norms with $p \in [1,\infty)$ and consider the absolute value instead of the positive part if we treat solutions.

The contributions of this article are twofold. On the one hand, hypoelliptic Poincar\'e inequalities are an interesting research topic and on the other hand, the emerging kinetic trajectories are interesting on their own. Let us comment on how these results and their proofs are connected to the hypoelliptic structure of the Kolmogorov equation. Consider $A = \id_d$. We write $X_0 = \partial_t + v \cdot \nabla_x $ and $X_i = \partial_{v_i}$, $i =1,\dots,d$. Then the Kolmogorov equation can be written as 
\begin{equation*}
	\sum_{i = 1}^d X_i^2 u - X_0u = 0,
\end{equation*}
i.e. as a H\"ormander sum of squares equation with drift $X_0$. The commutators
\begin{equation*}
 	[X_i,X_0]	=[\partial_{v_i},\partial_t + v \cdot \nabla_x ] = \partial_{v_i}(\partial_t + v \cdot \nabla_x ) - (\partial_t + v \cdot \nabla_x) \partial_{v_i} = \partial_{x_i},
\end{equation*}
$i = 1,\dots,d$, together with $X_0,X_1,\dots X_d$ generate a Lie-algebra of dimension $1+2d$ at every point $(t,x,v) \in \R^{1+2d}$. Thus the so-called H\"ormander's rank condition is satisfied, which implies hypoellipticity, \cite{hormander_hypoelliptic_1967}.

Let us recall some results on general H\"ormander vector fields $Y_0,Y_1,\dots,Y_m$ on $\R^N$ and the corresponding partial differential equation $\sum_{i = 1}^m Y_i^2u-Y_0u = 0$. Concerning trajectories along the vector fields connecting two points in space-time, it has already been observed in \cite{chow_systeme_1940} that if H\"ormander's rank condition is satisfied, then every two points can be connected by a trajectory along the vector fields $Y_0,Y_1,\dots,Y_m$. 

In the special case $Y_0 = 0$, the situation is well understood. In this case it is shown in \cite{nagel_balls_1985} that trajectories along the vector fields $Y_1,\dots,Y_m$ are equivalent to trajectories along the vector fields $Y_1,\dots,Y_m, [Y_i,Y_j],[Y_i,[Y_j,Y_k]], \dots $ and possibly higher order commutators. Equivalence has to be understood in the sense that they define equivalent metrics. The trajectories for $Y_0$ are very well-behaved and can be used to prove a Poincar\'e inequality, see \cite{lanconelli_poincare_2000}. Moreover, we refer to \cite{lu_horm_1992} and the references therein. If $Y_0 = \partial_t$, then the situation is very well understood and a Harnack inequality was proven in \cite{lu_horm_1992}. For general $Y_0$, there are fewer results of this spirit, and the situation is generally much more complicated. Recent progress towards a Harnack\rev{-type} inequality and a priori H\"older regularity for general H\"ormander sum of squares have been made in \cite{dietert_regularity_2022,citii_holder_2022,wang_hypo_2008}. 

We emphasise that there are some results, see \cite{pascucci_harnack_2004} which prove the existence of a \rev{kinetic} trajectory moving along the vector fields $X_0,X_1,\dots,X_d$ (corresponding to the Kolmogorov equation) connecting any two points. For the proof of Theorem \ref{thm:poincare}, we need a much more careful construction to ensure that the trajectory satisfies some beneficial properties. 
 
We want to compare the approach of \cite{guerand_quant_2021} with the one proposed in this article in view of its applicability to general H\"ormander sum of squares. Even if the vector fields $Y_0,Y_1,\dots,Y_d$ behave nicely, the commutators might be a mess and, hence, could be difficult to work with. However, the trajectories along the vector fields and their commutators could be of a simpler structure due to the increase in flexibility compared to the trajectories along $Y_0,Y_1,\dots,Y_d$. Nevertheless, results on the transfer of regularity are needed, which are currently not available in the general setting. As our proof does not depend on this kind of regularity transfer result, we hope that it can be extended to the general H\"ormander sum of squares. It is desirable to understand the trajectories in a much more general setting. An important question is for example to what extent the results of \cite{nagel_balls_1985} can be extended to $Y_0 \neq 0$.

We point out that we do not require any knowledge of the fundamental solution in our proof of Theorem \ref{thm:poincare} and use only trajectories corresponding to the geometry of the equation. Trajectorial arguments do not only allow relating the geometry of the equation to regularity properties but also allow for simpler proofs.  For instance, parabolic trajectories were studied in \cite{nz_moser_2022} to obtain a new interpretation of Moser's proof of the Harnack inequality. We hope that this trajectorial approach together with a better understanding of hypoelliptic trajectories allows proving the Harnack inequality for weak solutions to H\"ormander sum of squares with rough coefficients nicely. 

In Section \ref{sec:kkol}, we give some insight into how to transfer our method to the setting of H\"ormander sums of squares where commutators of higher order are needed to generate the Lie algebra. We investigate a Kolmogorov-type equation of order $k$ where commutators of order up to $k$ are needed to span the tangent space. \rev{We state the construction for arbitrary $k$, where we can verify the needed estimates for $k \le 9$.} To the knowledge of the authors, this is the first appearance of such a Poincar\'e inequality for hypoelliptic equations of order higher than two \rev{at the writing of the article}. \rev{We want to emphasise that this method has been successfully extended to arbitrary order of commutators in \cite{anceschi2024poincare} based on the construction proposed in this article.}

Concerning kinetic equations let us also note that in the work \cite{albritton_var_2019} a kinetic Poincar\'e inequality with an exponential weight in velocity is proven for solutions to the Kolmogorov equation with forcing term. This inequality, however, holds true only for solutions and it does not seem to be possible to use this inequality in the proof of the Harnack inequality. Moreover, it is not clear to the authors if their proof can be extended to hypoelliptic equations in general. Other Poincar\'e type inequalities in the kinetic setting appeared in the works \cite{anceschi_note_2022,guerand_log_2021,wang_nonhom_2007,wang_hypo_2008} where the function is compared not with a mean of the function but with a constant depending on the fundamental solution, too. Furthermore, trajectorial hypocoercivity has been used in \cite{dietert_trajectorial_2022} to study the long-time behaviour of solutions to kinetic equations. A kinetic Poincar\'e inequality is used in \cite{albritton_enh_2022} to prove enhanced dissipation. 

\rev{\textbf{Notation.} In what follows $d\in \N$ denotes the dimension. The subscripts $t$, $x$ and $v$ (of an $\R^{2d+1}$-valued function or of the gradient $\nabla$) refer to the first variable (time), the following $d$ (spatial) variables and the last $d$ (velocity) variables, respectively. Whenever we write $\partial_t +v\cdot \nabla_x$ we refer to the corresponding differential operator. We follow the standard notation for Lebesgue spaces $L^p$.}

\section{Kinetic trajectories}
\label{sec:curves}
In this section, we will construct \rev{the core object of this article, i.e.~kinetic} trajectories integral to the vector fields $\partial_t +v \cdot \nabla_x$ and $\partial_{v_i}$, $i = 1, \dots, d$. 

\rev{\begin{definition} \label{def:kintraj}
	Let $ (t,x,v),(s,y,w) \in \R^{1+2d}$ with $t \neq s$. A kinetic trajectory is a map $\gamma \in C([0,1];\R^{1+2d})$ differentiable in $(0,1)$ with $\gamma(0) = (t,x,v)$ and $\gamma(1) = (s,y,w)$ satisfying $\dot{\gamma}_x = \dot{\gamma}_t\gamma_v$ in $(0,1)$.
\end{definition}}

Let us first explain the idea of how to construct \rev{a kinetic trajectory $\gamma$ connecting two given points $ (t,x,v),(s,y,w) \in \R^{1+2d}$ with $s<t$.} For the time variable $\gamma_t$ \rev{the ansatz}
\begin{equation*}
	\gamma_t(r) = t+r^j(s-t)
\end{equation*}
with $j>0$ seems to be a good choice. We will later need that $\dot{\gamma}_t \le 0$ \rev{(as we are working with subsolutions)} so that, more generally, every decreasing function from $t$ to $s$ is possible, \rev{see also Remark \ref{rem:traj}}. 

Next we make \rev{an ansatz for the velocity component
\begin{equation} \label{eq:ansatzgammav}
	\gamma_v(r) = \alpha(r) c_1 + \beta(r)c_2+v
\end{equation}
for two continuous functions $\alpha,\beta\colon [0,1] \to \R$, differentiable in $(0,1)$} with the normalisation $\alpha(0) = \beta(0) = 0$ and $\alpha(1) = \beta(1) = 1$ and \rev{vectorial parameters} $c_1,c_2 \in \R^d$. 

\rev{The end point condition for the velocity component imposes} $w = \gamma(1) = c_1 +c_2 +v$ on the vectorial parameters. 

The condition that $\dot{\gamma}_x(r) = \dot{\gamma}_t(r)\gamma_v(r)$ (which makes sure that $\gamma$ only moves along $\partial_t+v \cdot \nabla_x$ and $\partial_{v_i}$) yields 
\begin{equation}
	\gamma_x(r) = \int_{0}^r \dot{\gamma}_t(\tau)\gamma_v(\tau) \dx \tau + x =  A(r) c_1 + B(r)c_2 + r^j(s-t)v + x
\end{equation}
where $A(r) = \int_{0}^r \dot{\gamma}_t(\tau) \alpha(\tau) \dx \tau $ and $B(r) = \int_{0}^r \dot{\gamma}_t(\tau) \beta(\tau) \dx \tau$ after integration in $r$. 

\rev{The endpoint condition $\gamma_x(1) = y$ imposes}
\begin{equation*}
	A(1)c_1 + B(1) c_2 + (s-t)v+x = y
\end{equation*}
as the second condition for $c_1,c_2$. If $A(1) \neq B(1)$ then we can explicitly solve \rev{the linear system
\begin{equation*}
	\begin{pmatrix}
		A(1) \; \id_d & B(1) \; \id_d \\
		\; \id_d & \; \id_d 
	\end{pmatrix} \begin{pmatrix}
		c_1 \\ c_2
	\end{pmatrix} = \begin{pmatrix}
		y-x-(s-t)v \\
		w-v
	\end{pmatrix}
\end{equation*}
}for $c_1$ and $c_2$ and obtain a kinetic trajectory in the sense of Definition \ref{def:kintraj}. \rev{As $c_1,c_2$ depend linearly on $x,v,y,w$ we may rewrite this kinetic trajectory as 
 \begin{equation} \label{eq:kintrajAB}
	\begin{pmatrix}
		\gamma_t(r) \\ \gamma_x(r)\\ \gamma_v(r)
	\end{pmatrix} = \begin{pmatrix}
		t+r^j(s-t)\\
		\A(r,s-t) \begin{pmatrix}
		y \\ w
	\end{pmatrix} + \B(r,s-t) \begin{pmatrix}
		x \\ v
	\end{pmatrix}
	\end{pmatrix}
\end{equation}}for two matrix-valued mappings $\A,\B \colon [0,1] \times \R \setminus \{0 \} \to \R^{2d \times 2d}$ satisfying $\A(0,s-t) = 0 = \B(1,s-t)$, $\A(1,s-t) = \id_{2d} = \B(0,s-t)$. We drop the $s-t$ dependence in our notation \rev{in the following}.

\rev{In general, these kinetic trajectories} are not good enough to make the proof of the Poincar\'e inequality work. We need to choose $j$ and $\alpha,\beta$ carefully. 

\rev{Let us first give our working choice and then explain how we derived it below in Remark \ref{rem:constrA}. We make the choice $j = 2$ based on the kinetic scaling ($\lambda \mapsto (\lambda^2 t, \lambda^3 x, \lambda v)$) which sets $\gamma_t(r) = t+(s-t)r^2$.} Then, we set
\begin{equation*}
	\A(r) = \begin{pmatrix}
		(7r^3-6r^{7/2}) \, \id_d & 4(\sqrt{r}-1)r^3(s-t) \, \id_d \\
		\frac{21}{2(s-t)}(r-r^{3/2}) \, \id_d & (7r^{3/2}-6r) \, \id_d \\
			\end{pmatrix}	
\end{equation*}
and
\begin{equation*}
	\B(r) = \begin{pmatrix}
	(1 + (6 \sqrt{r}-7) r^3) \, \id_d & (\sqrt{r}-1)^2 (1 + 2 \sqrt{r}) r^2 (s - t) \, \id_d \\
    \frac{21}{2(s-t)}(-r + r^{3/2}) \, \id_d & (1 - \frac{9r}{2} + \frac{7}{2}r^{3/2}) \, \id_d \\
			\end{pmatrix},
\end{equation*}
\rev{which together with equation \eqref{eq:kintrajAB} defines a kinetic trajectory.} It is not difficult to check the following important properties: 
\begin{enumerate}[(i)]
	\item $\dot{\gamma}_x(r) = \dot{\gamma}_t(r) \gamma_v(r) $ for all $r \in (0,1)$,
	\item $\det(\A(r)) = r^{\frac{9}{2}d}$, $r \in [0,1]$,
	\item $\abs{(\A^{-1})_{\cdot 2}} \le C_0(1+\abs{s-t})r^{-3/2}$ for all $r \in (0,1]$ some constant $C_0>0$,
	\item $\det(\B(r)) = \left(\frac{1}{2} (\sqrt{r}-1)^4 (1 + \sqrt{r}) (2 + 6 \sqrt{r} + 5 r + 6 r^{3/2} + 2 r^2)\right)^d$, for all $r\in [0,1]$,
	\item $\abs{\dot{\gamma}_v(r)} \le C_1\left(\frac{1}{\abs{s-t}}(\abs{x}+\abs{y})+\abs{v}+\abs{w}\right)$ for a universal constant $C_1>0$ and $r \in [0,1]$,
	\item and, if there is $\kappa>0$ such that $\kappa \le t-s \le 2+ \kappa$ and $x,v,y,w \in B_1(0) \subset \R^d$, then there exists a constant $R_0 \in [1,7]$ such that $(\gamma_x(r),\gamma_v(r)) \in B_{R_0(1+\kappa)}(0) \times B_{R_0(1+\kappa^{-1})}(0)$ for all $r \in [0,1]$. 
\end{enumerate}
\rev{We write $(\A^{-1})_{\cdot 2}$ for the $2d \times d$ matrix consisting of the upper right and lower right $d \times d$ block of $\A^{-1}$.  }

\begin{remark} \label{rem:constrA}
	The main difficulty is to find an appropriate $\A$ satisfying (iii) and (v) \rev{at the same time}. \rev{To find a kinetic trajectory with the good properties (ii)-(vi) we consider a reformulation of our above ansatz to simplify calculations. It is in fact equivalent to choosing $\alpha$ and $\beta$ in \eqref{eq:ansatzgammav}.} 
	
Let us assume that $\A$ is a $2d \times 2d$ block matrix of the form 
\begin{equation*}
	\A(r) = \begin{pmatrix}
		a(r) \, \id_d & b(r) \, \id_d \\
		c(r) \, \id_d & d(r) \, \id_d \\
	\end{pmatrix}
\end{equation*}
for real functions $a,b,c,d \colon [0,1] \to \R$. If additionally $\A(r)$ is assumed to be invertible for $r>0$, then 
\begin{equation*}
	\A^{-1}(r) = \frac{1}{a(r)d(r)-b(r)c(r)} \begin{pmatrix}
		d(r) \, \id_d & -b(r) \, \id_d \\
		-c(r) \, \id_d & a(r) \, \id_d \\
	\end{pmatrix}.
\end{equation*}
\rev{In particular, $(\A^{-1})_{\cdot 2} = (ad-bc)^{-1} (-b \, \id_d, a \, \id_d)^T$.} 

\rev{Let us explain the intuition behind the above construction.} The assumption $\dot{\gamma}_x(r) = \dot{\gamma}_t(r)\gamma_v(r)$ and coefficient comparison leads to $\dot{a}(r) = 2(s-t)rc(r)$ and $\dot{b}(r) = 2(s-t)rd(r)$. We need $c(0) = 0$ and $d(0) = 0$, hence $\dot{a}(0) = 0$ and $\dot{b}(0) = 0$. The function $c(r) = r^k$ for some $k>0$ seems to be a \rev{good ansatz}. Then, $a(r) = \int_0^r 2(s-t) \tau c(\tau) \dx \tau = \frac{2}{k+2} (s-t) r^{k+2}$. 

Recalling that $a$ describes the trajectory in the $x$-variable and $c$ the trajectory in the $v$-variable together with the \rev{kinetic} scaling \rev{($\lambda \mapsto (\lambda^2 t, \lambda^3 x, \lambda v)$)}, it seems natural to choose $k = 1$. This leads to $c(r) = r$ and $a(r) = 2/3(s-t)r^3$. The same \rev{reasoning} applies to $b$ and $d$. However, we cannot choose $a(r) = b(r)$ since otherwise the invertibility of $\A$ fails. 

Hence, we need to perturb $b$ a little bit but need to ensure that $(\A^{-1})_{\cdot 2}$ does not behave too badly. Thus, we made the ansatz $c(r) = c_1 r+c_2r^k$ and $d(r) = d_1 r+d_2r^k$ for some constants $c_1,c_2,d_1,d_2 \in \R$. The parameters $c_1,c_2,d_1,d_2$ are determined by the \rev{start} and \rev{end point} conditions. \rev{We need $k \ge 1$ in order for (v) to hold and $k \neq 1$ in light of the above argument.} Choosing $k = 3/2$ leads to $\A$ as above. \rev{Any choice of} $k >1$ is possible and would lead to \rev{the determinant scale like $r^{(3+k)d}$} and a different decay in the property (iii) of order \rev{$r^{-k}$. In the proof of Theorem \ref{thm:poincare} we need that $r^{1-k}$ is integrable close to zero which is why we need to restrict to $k \in (1,2)$.} 

The matrix $\B$ can be constructed similarly. We note that infinitely many choices exist for the mappings $\A$ and $\B$. By tweaking the matrices $\A$ and $\B$, one can try to optimise the constants $C_0,C_1,R_0$.
\end{remark}

\begin{remark} \label{rem:traj} 
$ $ 
	\begin{enumerate}[1)]
		\item 	Note that with a suitable $\A$ at hand and for $\B= \B(r) \in C^1$ on $[0,1]$ satisfying (i), the properties (v),(vi) are immediate, and $\det(\B(r))\ge \tilde{C}>0$ on $[0,\eta]$ for a possibly small $\eta>0$ (this is a weaker version of (iv) which suffices for our approach), since $\B(0) = \id$ and by continuity. To be more precise, in neither (ii) nor (iv) the precise formula is important for later calculations but only the positivity in suitable intervals is needed.  
		\item One could choose a different ansatz for $\gamma_t(r)$ (for example $j \neq 2$), which in turn leads to different $\A(r)$ and $\B(r)$.
		\item In view of \cite{nz_moser_2022} it is an interesting open problem whether one can construct a trajectory mapping $\A$ where (iii) can be improved to $\abs{(\A^{-1})_{\cdot 2}} \le C_0r^{-1}$. However, for smooth $a,b,c,d$ one can prove by a Taylor expansion argument that the trajectory cannot reach this behaviour. In the parabolic case, it is easy to construct a trajectory with this behaviour, compare \cite{nz_moser_2022}. Understanding this phenomenon better is an important step towards a better understanding of the hypoelliptic De Giorgi-Nash-Moser theory. 
		\item Apparently the constructed trajectories look rather complicated. It would be desirable to have a better understanding of them for example by an abstract principle or by understanding the exponential map of the vector fields. 
	\end{enumerate}
\end{remark}

\section{Proof of Theorem \ref{thm:poincare}}

\begin{proof}
	Recall that $Q_1^- =(-2-\kappa,-1-\kappa] \times B_1 \times B_1$ and abbreviate $B = B_1 \times B_1$, where $B_1 = B_1(0)$. Let $\varphi \in C_c^\infty(B_1 \times B_1)$ such that $\varphi = 1$ in $B_{1-\epsilon} \times B_{1-\epsilon}$, $0 \le \varphi \le 1$ and $\abs{\nabla \varphi} \le C_2 \epsilon^{-1}$ for some $C_2>0$. 
	
To start we proceed as in \cite{guerand_quant_2021} and estimate
	\begin{align*}
		&\norm{(u-\langle u \rangle_{Q_1^-})_+}_{L^1(Q_1)} \le \norm{(u - \langle u \varphi^2 \rangle_{Q_1^-})_+}_{L^1(Q_1)} \\
		 &\le \int_{Q_1} \left( \fint_{Q_1^{-}} (u(t,x,v)-u(s,y,w)) \varphi^2(y,w) \dx (s,y,w) \right)_+ \dx (t,x,v) \\
		 &\hphantom{=} +  \norm{u}_{L^1(Q_1)} \fint_{Q_1^-}(1-\varphi^2(y,w)) \dx (s,y,w) \\
		 &\le \int_{Q_1} \left( \fint_{Q_1^{-}} (u(t,x,v)-u(s,y,w)) \varphi^2(y,w) \dx (s,y,w) \right)_+ \dx (t,x,v) + C\epsilon^{2d} \norm{u}_{L^1(Q_1)}  
	\end{align*}
	where $C = C(d)$. 

We write the integrand of the first term as follows
\begin{align} \label{eq:estimate}
	&\int_{Q_1^{-}} (u(t,x,v)-u(s,y,w)) \varphi^2(y,w) \dx (s,y,w) \\
	&= - \int_{Q_1^{-} }\int_0^1 \frac{\dx }{\dx r} u (\gamma(r)) \dx  r \varphi^2(y,w) \dx (s,y,w) \nonumber \\
	&= -  \int_{Q_1^{-}} \int_0^1 \dot{\gamma}_t(r)[(\partial_t + v \cdot \nabla_x) u](\gamma(r)) \dx r \varphi^2(y,w) \dx (s,y,w) \nonumber \\
	&\hphantom{=} -   \int_{Q_1^{-}}  \int_0^1 \dot{\gamma}_v(r)\cdot[\nabla_v u](\gamma(r)) \dx r \varphi^2(y,w) \dx (s,y,w) =: I_1 + I_2, \nonumber
\end{align} 	
where $\gamma \colon [0,1] \to \R^{1+2d}$ is a trajectory going from $(t,x,v)$ to $(s,y,w)$ as constructed in Section \ref{sec:curves}. Note that we used property (i) in the last equality. 

We continue by estimating the first term. Recall that $\dot{\gamma}_t(r) = 2r(s-t) \le 0$ as $s-t<0$. Hence, we can use the subsolution property of $u$. Next, apply Fubini's theorem and  substitute $(\tilde{y},\tilde{w}) = \Phi(y,w)=\A(r) (y,w)^T + \B(r)(x,v)^T$, where we think of $s,t,x,v,r$ as fixed parameters and $(y,w)$ \rev{as variables}. We obtain
\begin{align*}
	I_1 &\le -\int_0^1\int_{Q_1^{-}}  \dot{\gamma}_t(r)[\nabla_v \cdot (A\nabla_v u)](\gamma(r))\varphi^2(y,w)  \dx (s,y,w)\dx r \\
	&=-\int_{-2-\kappa}^{-1-\kappa} \int_0^1\int_{\Phi^{-1}(B)} \dot{\gamma}_t(r)[\nabla_v \cdot (A\nabla_v u)](\gamma_t(r),\tilde{y},\tilde{w}) \varphi^2(\Phi^{-1}(\tilde{y},\tilde{w})) \\
	&\hspace{9cm} \abs{\det{\A(r)}}^{-1} \dx (\tilde{y},\tilde{w}) \dx r  \dx s \\
	&= \int_{-2-\kappa}^{-1-\kappa}2\int_0^1 \int_{\Phi^{-1}(B)}  \dot{\gamma}_t(r)[A\nabla_v u](\gamma_t(r),\tilde{y},\tilde{w}) \cdot (\nabla \varphi(\Phi^{-1}(\tilde{y},\tilde{w})) (\A^{-1})_{\cdot 2} )  \varphi(\Phi^{-1}(\tilde{y},\tilde{w})) \\
	&\hspace{10.5cm} \abs{\det{\A(r)}}^{-1} \dx (\tilde{y},\tilde{w}) \dx r  \dx s \\
	&= 2\int_{Q_1^-} \int_0^1 \dot{\gamma}_t(r)[A\nabla_v u](\gamma(r)) \cdot (\nabla \varphi(y,w) (\A^{-1})_{\cdot 2} )  \varphi(y,w) \dx r \dx (s,y,w)  \\
	&\le \Lambda \epsilon^{-1} C_3 \int_{Q_1^-} \int_0^1 r^{-1/2}\abs{\nabla_v u}(\gamma(r))  \dx r \dx (s,y,w) 
\end{align*}
with $C_3 = C_3(d)>0$ by the properties (ii) and (iii) written in Section \ref{sec:curves}. In the third line, we performed a partial integration, recalling that $\Phi$ is a homeomorphism $\Phi (\partial B) = \partial \Phi(B)$ so that the boundary term vanishes.
A straightforward estimate of the second term yields
\begin{align*}
	I_2 &\le \int_{Q_1^-}  \int_0^1 \abs{\dot{\gamma}_v(r)}\abs{\nabla_v u}(\gamma(r)) \dx r \varphi^2(y,w) \dx (s,y,w) \\
	&\le C_4\int_{Q_1^-}  \int_0^1 \abs{\nabla_v u}(\gamma(r)) \dx r \dx (s,y,w)
\end{align*}
for some constant $C_4 = C_4(\kappa)>0$ as a consequence of property (v) from Section \ref{sec:curves}.

Let $k>-1$. We want to estimate the integral
\begin{align*}
	J&:=\int_{Q_1} \int_{-2-\kappa}^{-1-\kappa}\int_{B} \int_0^1 r^{k}\abs{\nabla_v u}(\gamma(r))  \dx r \dx (y,w) \dx s \dx (t,x,v) \\
	&= \int_{Q_1} \int_{-2-\kappa}^{-1-\kappa}\int_{B} \int_0^{1/2} r^{k}\abs{\nabla_v u}(\gamma(r))  \dx r \dx (y,w) \dx s \dx (t,x,v)  \\
	&\hphantom{=}+ \int_{Q_1} \int_{-2-\kappa}^{-1-\kappa}\int_{B} \int_{1/2}^1 r^{k}\abs{\nabla_v u}(\gamma(r)) \dx r \dx (y,w) \dx s \dx (t,x,v) =: J_1+J_2,
\end{align*}
where we split the $r$-integral on the interval $[0,1]$ into the integrals on $[0,1/2]$ and $[1/2,1]$. 

For the first integral $J_1$ we apply Fubini's theorem, substitute $\tilde{t} = t+r^2(s-t)$ and $(\tilde{x},\tilde{v}) = \Psi(x,v)=\A(r) (y,w)^T + \B(r)(x,v)^T$ and estimate as follows
\begin{align*}
	J_1 &= \int_{-2-\kappa}^{-1-\kappa} \int_{B} \int_0^{1/2} \int_{-1}^0 \int_{B} r^{k}\abs{\nabla_v u}(\gamma(r))  \dx (t,x,v) \dx r \dx (y,w) \dx s  \\
	&= \int_{-2-\kappa}^{-1-\kappa} \int_{B} \int_0^{1/2} \int_{r^2s+(r^2-1)}^{r^2s} \int_{\Psi(B)} r^{k}\abs{\nabla_v u}(\tilde{t},\tilde{x},\tilde{v}) \abs{\det(\B)}^{-1} \\
	&\hphantom{=} \hspace{6.5cm} (1-r^2)^{-1} \dx (\tilde{x},\tilde{v}) \dx \tilde{t}  \dx r \dx (y,w) \dx s  \\
	&\le C_5  \int_{-2-\kappa}^{-1-\kappa} \int_{B} \int_0^{1/2} \int_{-2-\kappa}^0 \int_{\tilde{B}} r^{k}\abs{\nabla_v u}(\tilde{t},\tilde{x},\tilde{v}) \dx (\tilde{x},\tilde{v}) \dx \tilde{t}  \dx r \dx (y,w) \dx s\\
	&\le C_5 \int_{-2-\kappa}^{-1-\kappa} \int_{B}  \int_{\tilde{Q}}  \int_0^{1/2} r^{k} \dx r \abs{\nabla_v u}(\tilde{t},\tilde{x},\tilde{v}) \dx (\tilde{t},\tilde{x},\tilde{v})  \dx (y,w) \dx s \le \abs{Q_1} C_6 C_5 \norm{\nabla_v u}_{L^1(\tilde{Q})}.
\end{align*}
where $\abs{\det (\B)}^{-1} \le C_5$ on $[0,1/2]$ due to property (iv) in Section \ref{sec:curves}, $\tilde{B} = B_{R_0}(0) \times B_{R_0(1+\kappa^{-1})}(0)$ with $\Psi(B) \subset \tilde{B}$ by property (vi) as well as $\tilde{Q} = (-2-\kappa,0) \times \tilde{B}$ and $C_6=C_6(k)$.

For the term $J_2$ we substitute $\tilde{s} = t+r^2(s-t)$ and $(\tilde{y},\tilde{w}) = \Phi(y,w)=\A(r) (y,w)^T + \B(r)(x,v)^T$ to obtain
\begin{align*}
	J_2 &= \int_{Q_1} \int_{1/2}^1 \int_{-2-\kappa}^{-1-\kappa} \int_{B} r^{k}\abs{\nabla_v u}(\gamma(r)) \dx (y,w) \dx s \dx r \dx (t,x,v)  \\
	&= \int_{Q_1} \int_{1/2}^1 \int_{(1-r^2)t+r^2(-2-\kappa)}^{(1-r^2)t+r^2(-1-\kappa)} \int_{\Phi(B)}  r^{k-2}\abs{\nabla_v u}(\tilde{s},\tilde{y},\tilde{w}) \dx (\tilde{y},\tilde{w})  \dx \tilde{s}  \abs{\det(\A)}^{-1} \dx r  \dx (t,x,v) \\
	&\le \int_{Q_1} \int_{1/2}^1 \int_{-2-\kappa}^0 \int_{\tilde{B}}  r^{k-2-(9/2)d}\abs{\nabla_v u}(\tilde{s},\tilde{y},\tilde{w}) \dx (\tilde{y},\tilde{w})  \dx \tilde{s}   \dx r  \dx (t,x,v) \le C_7 \abs{Q_1} \norm{\nabla_v u}_{L^1(\tilde{Q})},
	\end{align*}
	by property (ii) in Section \ref{sec:curves}, where $C_7 = C_7(d,k)$.
	
	Using the estimate for $J$ with $k = -1/2$ and $k = 0$ we conclude
	\begin{align*}
		&\int_{Q_1} \left( \fint_{Q_1^{-}} (u(t,x,v)-u(s,y,w)) \varphi^2(y,w) \dx (s,y,w) \right)_+ \dx (t,x,v) \\
		&\hphantom{=}\le \abs{Q_1^{-}}^{-1} \int_{Q_1} I_1 \dx (t,x,v) + \abs{Q_1^{-}}^{-1} \int_{Q_1} I_2 \dx (t,x,v)  \\
		&\hphantom{=}\le C\epsilon^{-1} \norm{\nabla_v u}_{L^1(\tilde{Q})},
	\end{align*}
	where $C = C(d,\Lambda,\kappa)$ and $\epsilon \in (0,1)$. This concludes the proof of the theorem. 
	
	As in \cite{guerand_quant_2021}, we use formal calculations. For instance, the estimate for the term in equation \eqref{eq:estimate} needs to be understood in a weak sense, i.e. the trajectorial argument needs to be performed on the level of test functions.
\end{proof}

\begin{remark}
	Looking at the proof it is immediately clear that this trajectorial proof of the Poincar\'e inequality holds true also for more general hypoelliptic equations if suitable trajectories are available. The difficulty lies in constructing trajectories enjoying properties similar to (i)-(vi) of Section \ref{sec:curves}. To the knowledge of the authors, this is the first time this argument has been successfully employed in the hypoelliptic setting with drift. 
\end{remark}

\section{Towards more general hypoelliptic equations}
\label{sec:kkol}
In this section we consider solutions $u = u(t,x_1,\dots,x_k,v) \colon \R^{1+kd+d}$ for $k \ge 1$ to the equation 
\begin{equation} \label{eq:kolk}
	(\partial_t  +  \sum_{j = 2}^{k-1} x_{j} \nabla_{x_{j-1}} + v \cdot \nabla_{x_k}) u = \Delta_v u.
\end{equation}
This equation is hypoelliptic and serves as a natural generalisation of the Kolmogorov equation. This time the commutators of orders up to $k$ are needed to generate the Lie algebra of dimension $1+(k+1)d$. The equation is an important example to investigate whether our proposed method can be extended to more general hypoelliptic equations. We refer to \cite{anceschi_note_2022} and the references therein. 

We conjecture the following.
\begin{conj}\label{thm:poincare2}
	There exists a universal constant $R_0>0$ such that the following holds. Let $\epsilon \in (0,1)$, $\kappa>0$ and $u$ be a nonnegative weak subsolution to equation \eqref{eq:kolk} in $\tilde{Q} = (-2-\kappa,0] \times B_{R_0} \times B_{R_0(1+\kappa^{-1})} \dots  \times B_{R_0(1+\kappa^{-1})^{k}}$. Then 
	\begin{equation*}
		\norm{(u-\langle u \rangle_{Q_1^-})_+}_{L^1(Q_1)} \le C \epsilon^{-1} \norm{\nabla_v u}_{L^1(\tilde{Q})} + C\epsilon^{(k+1)d} \norm{u}_{L^1(Q_1)}
	\end{equation*}
	where $Q_1^- = (-2-\kappa,-1-\kappa] \times B_1 \times \dots  \times B_1$, $Q_1 = (-1,0] \times B_1 \times \dots  \times B_1$, $\langle u\rangle_{Q_1^-} = \fint_{Q_1^-} u \dx(t,x,v)  = \frac{1}{\abs{Q_1^-}} \int_{Q_1^-} u \dx (t,x,v)$, where $x = (x_1,\dots,x_k)$ and $C = C(d,\kappa)$.
\end{conj}

\begin{remark}
	\rev{After the release of this work as a preprint but prior to its publication, our construction of suitable hypoelliptic trajectories was successfully extended to any order of commutators in the preprint \cite{anceschi2024poincare}, where also a proof of Conjecture \ref{thm:poincare2} was given.}
\end{remark}

Let $(t,x,v) = (t,x_1,\dots,x_k,v) \in Q_1,(s,y,w) = (s,y_1,\dots,y_k,w) \in Q_1^-$. We are looking for a trajectory $\gamma = (\gamma_t,\gamma_{x_1},\dots,\gamma_{x_k},\gamma_v) \colon [0,1] \to \R^{1+(k+1)d}$ such that $\gamma(0) = (t,x,v)$ and $\gamma(1) = (s,y,w)$ with
\begin{equation} \label{eq:pathk}
	\dot{\gamma}_{x_k}(r) = \dot{\gamma}_t(r)\gamma_v(r), \quad \dot{\gamma}_{x_{j-1}}(r) = \dot{\gamma}_t(r)\gamma_{x_{j}}(r), \; j = 2,\dots,k   
\end{equation}
for all $r \in [0,1]$. This implies
\begin{equation*}
	\frac{\dx}{\dx r} u(\gamma(r)) = \dot{\gamma}_t(r)\left[(\partial_t + \sum_{j = 2}^{k-1} x_{j} \nabla_{x_{j-1}}  + v \cdot \nabla_{x_k}) u\right](\gamma(r)) + \dot{\gamma}_v(r) \cdot [\nabla_v u](\gamma(r))
\end{equation*}
for all sufficiently smooth $u \colon \R^{1+(k+1)d} \to \R$. If we are able to construct a trajectory with properties similar to (i)-(vi) in Section \ref{sec:curves}, then it is immediate that the proof of Theorem \ref{thm:poincare} carries over to prove Conjecture \ref{thm:poincare2}. We emphasise that it is easy to find a trajectory connecting two points with property (i), compare \cite{pascucci_harnack_2004}. The intricate part is finding a trajectory which satisfies properties analogous to (ii)-(vi).

Following the same ansatz as in Section \ref{sec:curves} we propose to construct a suitable trajectory as follows. Set $\gamma_t(r) = t+(s-t)r^2$ and $(\gamma_x(r),\gamma_v(r))^T = \A(r,s-t) (y,w)^T + \B(r,s-t)(x,v)^T$, with two matrix-valued mappings $\A,\B \colon [0,1] \times \R \setminus \{0 \} \to \R^{(k+1)d \times (k+1)d}$ satisfying $\A(0,s-t) = 0 = \B(1,s-t)$, $\A(1,s-t) = \id = \B(0,s-t)$. Recall that the main difficulty is to find a suitable mapping $\A$, see the comments in Section 2. 

Think of $\A$ as a block matrix of $(k+1)\times (k+1)$ identity matrices of dimension $d$ with scalar prefactors, i.e. $\A = (a_{ij}\id_d)_{i,j = 1,\dots,k+1}$. Set $(\kappa_{i})_{i = 1,\dots, k} = (1+\frac{1}{i+1})_{i = 1,\dots, k}$ and define
\begin{equation*}
	a_{(k+1)j} = \alpha_{1j} \ r+ \sum_{l = 2}^{k+1} \alpha_{lj} \ r^{\kappa_{l-1}}, \; j=1,\dots,k+1
\end{equation*}
for some constants $(\alpha_{ij})_{i,j = 1,\dots,k+1}$. This is consistent with our choice of $c$ and $d$ in Section \ref{sec:curves}. The relations in equation \eqref{eq:pathk} dictate the structure of $a_{ij}$ for $i = 1, \dots, k$ and $j = 1, \dots, k+1$ by iteratively multiplying $a_{(k+1)j}$ with $\dot{\gamma}_t$ and integrating from $0$ to $r$. Note that the initial condition $\A(0) = 0$ is satisfied by construction. The terminal condition $\A(1) = \id$ leads to a system of $(k+1)^2$-equations for the coefficients $(\alpha_{ij})_{i,j = 1,\dots,k+1}$. $\B$ is constructed similarly. We conjecture that $\A,\B$ satisfy comparable properties to those listed in Section \ref{sec:curves}, which would allow transferring the proof of Theorem \ref{thm:poincare} to Conjecture \ref{thm:poincare2}. 

With this ansatz, the problem boils down to a purely algebraic problem. We give evidence for $k =2 $ and $k = 3$ below. We pose the proof for general $k $ as an interesting open problem. Due to the structure of the mappings, properties (i),(ii),(iv),(v) and (vi) should be possible to check. The main difficulty lies in proving a property similar to (iii). We conjecture that $\abs{(\A^{-1})_{\cdot k+1}} \le Cr^{-3/2}$ for all $k \in \N$ and some $C>0$ and $\det(\A ) = r^{p_k d}$ with $p_k = (k+1)^2+ \sum_{l = 1}^k [(\kappa_l-1)] $. This holds true for $k = 2,3$ and has been verified for $k = 4,5,6,7,8,9$ with the help of computer algebra. 

We emphasise that this is just a particular choice for the trajectories and that there is a lot of flexibility concerning the temporal trajectory and the exponents $\kappa$. 

For $k =2$ and $k = 3$, the coefficients $(\alpha_{ij})_{i,j = 1,\dots,k+1}$ are uniquely determined, and analogous properties to (i) to (vi) of Section \ref{sec:curves} can be proven. In particular, $\abs{(\A^{-1})_{\cdot 3}} \le Cr^{-3/2}$ and $\abs{(\A^{-1})_{\cdot 4}} \le Cr^{-3/2}$, respectively hold for some constant $C>0$ and all $r \in (0,1]$. 

\textbf{The case $k = 2$.}
We obtain {\tiny 
\begin{equation*}
	\A(r) = \left(
\begin{array}{ccc}
 \left(-495 \sqrt[3]{r}+320 \sqrt{r}+176\right) r^5 & -2 \left(-153 \sqrt[3]{r}+100 \sqrt{r}+53\right) r^5 (s-t) & 24 \left(2 r^{11/2}-3 r^{16/3}+r^5\right) (s-t)^2 \\
 \frac{440 \left(2 r^{7/2}-3 r^{10/3}+r^3\right)}{s-t} & \left(816 \sqrt[3]{r}-550 \sqrt{r}-265\right) r^3 & 12 \left(-16 \sqrt[3]{r}+11 \sqrt{r}+5\right) r^3 (s-t) \\
 \frac{220 \left(-10 \sqrt[3]{r}+7 \sqrt{r}+3\right) r}{(s-t)^2} & -\frac{5 \left(-544 \sqrt[3]{r}+385 \sqrt{r}+159\right) r}{2 (s-t)} & \left(-320 \sqrt[3]{r}+231 \sqrt{r}+90\right) r \\
\end{array}
\right)
\end{equation*}}
and $\abs{(\A^{-1})_{\cdot 3}} \le Cr^{-3/2}$ for all $r \in (0,1]$ and some $C>0$ where $(\A^{-1})_{\cdot 3}$ calculates as{\small
\begin{equation*}
	\left(\frac{24 \left(-3 \sqrt[6]{r}+\sqrt{r}+2\right) (s-t)^2}{r^{3/2}},\frac{12 \left(-16 \sqrt[6]{r}+5 \sqrt{r}+11\right) (s-t)}{r^{3/2}},\frac{-320 \sqrt[6]{r}+90 \sqrt{r}+231}{r^{3/2}}\right)^T.
\end{equation*}}

\textbf{The case $k = 3$.}
We find {\scriptsize 
\begin{align*}
	\A(r) &= \left(
\begin{array}{cc}
 \left(82215 \sqrt[3]{r}-73920 \sqrt[4]{r}-17864 \sqrt{r}+9570\right) r^7 & 2 \left(-25515 \sqrt[3]{r}+22880 \sqrt[4]{r}+5572 \sqrt{r}-2937\right) r^7 (s-t)  \\
 -\frac{33495 \left(-9 \sqrt[3]{r}+8 \sqrt[4]{r}+2 \sqrt{r}-1\right) r^5}{s-t} & \left(-187110 \sqrt[3]{r}+165880 \sqrt[4]{r}+41790 \sqrt{r}-20559\right) r^5  \\
 -\frac{33495 \left(-48 \sqrt[3]{r}+42 \sqrt[4]{r}+11 \sqrt{r}-5\right) r^3}{2 (s-t)^2} & \frac{1155 \left(-864 \sqrt[3]{r}+754 \sqrt[4]{r}+199 \sqrt{r}-89\right) r^3}{2 (s-t)} \\
 -\frac{33495 \left(-320 \sqrt[3]{r}+273 \sqrt[4]{r}+77 \sqrt{r}-30\right) r}{8 (s-t)^3} & \frac{1155 \left(-5760 \sqrt[3]{r}+4901 \sqrt[4]{r}+1393 \sqrt{r}-534\right) r}{8 (s-t)^2} \\
\end{array} 
\right. \\
&\left.\begin{array}{cc}
 -8 \left(-1701 \sqrt[3]{r}+1520 \sqrt[4]{r}+374 \sqrt{r}-193\right) r^7 (s-t)^2 & 192 \left(-9 \sqrt[3]{r}+8 \sqrt[4]{r}+2 \sqrt{r}-1\right) r^7 (s-t)^3 \\
  -4 \left(-12474 \sqrt[3]{r}+11020 \sqrt[4]{r}+2805 \sqrt{r}-1351\right) r^5 (s-t) & 96 \left(-66 \sqrt[3]{r}+58 \sqrt[4]{r}+15 \sqrt{r}-7\right) r^5 (s-t)^2 \\
 \left(133056 \sqrt[3]{r}-115710 \sqrt[4]{r}-30855 \sqrt{r}+13510\right) r^3 & 24 \left(-704 \sqrt[3]{r}+609 \sqrt[4]{r}+165 \sqrt{r}-70\right) r^3 (s-t) \\
 -\frac{105 \left(-8448 \sqrt[3]{r}+7163 \sqrt[4]{r}+2057 \sqrt{r}-772\right) r}{4 (s-t)} & \left(-28160 \sqrt[3]{r}+23751 \sqrt[4]{r}+6930 \sqrt{r}-2520\right) r \\
\end{array} \right) \end{align*}}
and $\abs{(\A^{-1})_{\cdot 4}} \le Cr^{-3/2}$ for all $r \in (0,1]$ and some $C>0$, where $(\A^{-1})_{\cdot 4}$ calculates as{ \small 
\begin{align*} 
	&\left(-\frac{192 \left(-8 \sqrt[4]{r}+9 \sqrt[6]{r}+\sqrt{r}-2\right) (s-t)^3}{r^{3/2}},-\frac{96 \left(-58 \sqrt[4]{r}+66 \sqrt[6]{r}+7 \sqrt{r}-15\right) (s-t)^2}{r^{3/2}},\right. \\
	&\hphantom{=}\left.-\frac{24 \left(-609 \sqrt[4]{r}+704 \sqrt[6]{r}+70 \sqrt{r}-165\right) (s-t)}{r^{3/2}},\frac{23751 \sqrt[4]{r}-28160 \sqrt[6]{r}-2520 \sqrt{r}+6930}{r^{3/2}}\right)^T.
\end{align*}}

\noindent \rev{\textbf{Acknowledgements.} The authors wish to thank Cyril Imbert for helpful discussions and hosting us at \'Ecole normale sup\'erieure in November 2021 where part of this research was initiated. Moreover, we thank Helge Dietert and Cl\'ement Mouhot for fruitful discussions. Lukas Niebel is funded by the Deutsche Forschungsgemeinschaft (DFG, German Research Foundation) under Germany's Excellence Strategy EXC 2044 --390685587, Mathematics M\"unster: Dynamics--Geometry--Structure.}

\bibliographystyle{amsplain}

\providecommand{\bysame}{\leavevmode\hbox to3em{\hrulefill}\thinspace}
\providecommand{\MR}{\relax\ifhmode\unskip\space\fi MR }

\providecommand{\MRhref}[2]{%
  \href{http://www.ams.org/mathscinet-getitem?mr=#1}{#2}
}
\providecommand{\href}[2]{#2}

$ $\\

\end{document}